\documentclass{amsart}
\usepackage{graphicx}
\usepackage[pdfauthor={Richard J. Mathar},pdfkeywords={Gauss Quadrature, Analysis, Logarithmic Kernel},pdfsubject={Mathematics, Analysis},colorlinks=true,citecolor=blue]{hyperref}

\usepackage{amsmath}

\newtheorem{rem}{Remark}
\newtheorem{exa}{Example}
\def\Cl{\mathop{\mathrm{Cl}}\nolimits}

\begin{document}

\title[C Implementation of Clausen sums]{A C99 Implementation of the Clausen Sums.}

\author{Richard J. Mathar} 
\email{mathar@mpia.de}
\urladdr{http://www.mpia.de/~mathar}
\address{Hoeschstr. 7, 52372 Kreuzau, Germany}

\subjclass[2010]{Primary 33E20, 33F05; Secondary 26-04, 65D20}

\date{\today}
\keywords{Clausen Integral, Special Functions, Log-Sine Integral}

\begin{abstract}
A C99 program is presented which calculates the Clausen sums for non-negative
integer indices
$j$ and real arguments $x$. The implementation is split into three major cases.
For large $j$, a direct, truncated summation of the sum is performed. For smaller $j$,
the associate polynomial of $x$ is evaluated if applicable; otherwise a 
Chebyshev series is constructed by repeated integration of a
(truncated) Chebyshev representation
of $\Cl_2(x)$, and evaluated.
\end{abstract}

\maketitle 

\section{Definitions} 

The theme of this paper is the numerical evaluation of the two 
classes of sums defined as
\begin{equation}
C_j(x)\equiv \sum_{k\ge 1}\frac{\cos(kx)}{k^j}
\label{eq.Cdef}
\end{equation}
and
\begin{equation}
S_j(x)\equiv \sum_{k\ge 1}\frac{\sin(kx)}{k^j}
\label{eq.Sdef}
\end{equation}
for integers $j\ge 1$ and arguments $x$ on the real line.
The standard definition of  Clausen sums $\Cl_j(x)$
addresses the sine type for even $j$
and the cosine type for odd $j$:
\begin{equation}
\Cl_{2j}(x)\equiv S_{2j}(x);\quad
\Cl_{2j+1}(x)\equiv C_{2j+1}(x).
\label{eq.Clausdef}
\end{equation}

\begin{figure}
\includegraphics[width=0.8\columnwidth]{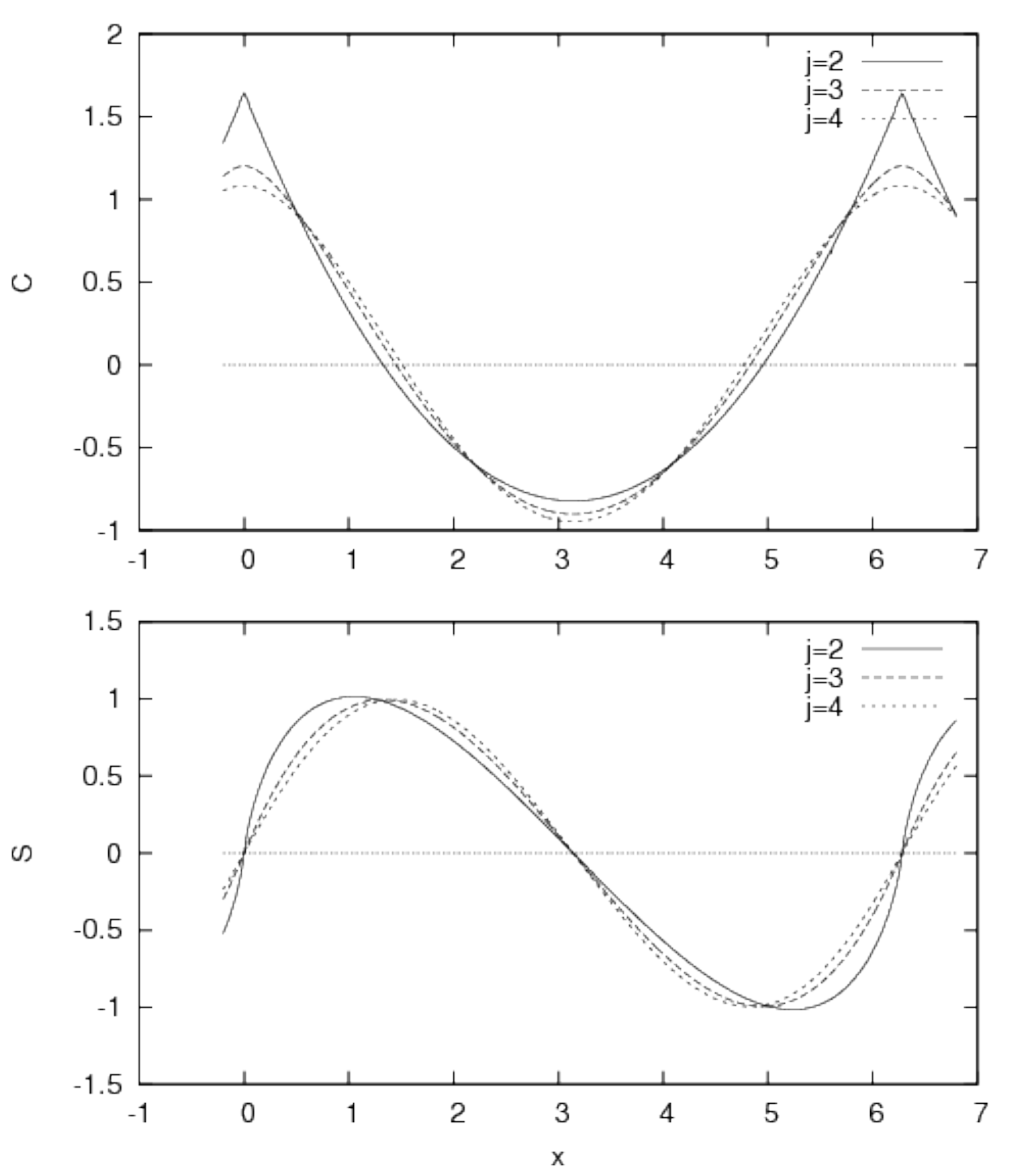}
\caption{
The functions $C_j(x)$ (upper window) and $S_j(x)$ (lower window) for small $j$.
The display covers a little bit more than one $2\pi$ period.
}
\label{fig.CSofj}
\end{figure}

The values are periodic in $x$ with a period of $2\pi$,
illustrated in Figure \ref{fig.CSofj}\@.
$C_2(x)$ is continuous but not differentiable where $x$
is a multiple of $2\pi$.
$S_1(x)$ is a saw-tooth function, discontinuous at places where $x$
is a multiple of $2\pi$.

Special values have been discussed in the literature if
$x$ is a rational multiple of $\pi$
\cite{CvijovicMathComp64,DolderJCAM11,YangIJMEST23,SrivastavaZAA19}.
This did no take impact on the current implementation.

\section{Representation and Approximation}
\subsection{Large index or index equal to 1}
A simple implementation is available for large $j$, which is actually
implemented for $j>10$.
Since the sums are converging quickly as a result of the high powers,
a partial sum of some dozen of the first terms suffices for
double precision accuracy.

For $j=1$, the closed forms
\begin{equation}
C_1(x)= -\log\left(2\sin\frac{x}{2}\right),\quad 0<x<2\pi,
\label{eq.C1}
\end{equation}
\begin{equation}
S_1(x)= \frac{1}{2}(\pi-x),\quad 0<x<2\pi
\label{eq.S1}
\end{equation}
are used \cite[27.8.6]{AS}\cite{ClausenCrelle8}\cite[4.1.3]{Apelblat2}.

The remaining cases of small $j$ split into two classes which are discussed
in the following two subsections, first the simpler combinations of sine/cosine
and parity of $j$ in Section \ref{sec.poly},
then the combinations related to the functions (\ref{eq.Clausdef})
in Section \ref{sec.Cl}.

\subsection{Polynomial Branch} \label{sec.poly}
The functions $S_j(x)$ are polynomials of $x$
if $j$ is odd, and the functions $C_j(x)$ are polynomials
of $x$ if $j$ is even \cite{Jolley}\cite[27.8.6]{AS}. The 
basic relation is
\begin{equation} 
C_2(x)=\frac{\pi^2}{6}-\frac{\pi}{2}x+\frac{1}{4}x^2,\quad 0\le x\le 2\pi.
\label{eq.C2pol}
\end{equation} 
By successive integration with respect to $x$ with lower limit at $x=0$, with interchange
of summation and integration, the
functions $C_2(x)\to S_3(x) \to C_4(x) \to S_5(x)\to \ldots$ are also
represented as polynomials of order $j$ \cite{BoydJCP228}. The constant of the integration
is fixed by enforcing the limits $S_j(0)=0$ and $C_j(0)=\zeta(j)$.
\begin{rem}
The Riemann Zeta Function $\zeta(j)$ is a fractional multiple of $\pi^j$ in
these cases.
\end{rem}
\begin{rem}
For each step $S_j(x) \to C_{j+1}(x)$ an additional
sign flip is induced, because
$\int \sin (kx)/k^j dx= - \cos(kx)/k^{j+1}+\mathrm{const}$.
\end{rem}
\begin{exa}
This yields \cite[27.8.6]{AS}\cite[p. 96]{Jolley}
\begin{eqnarray} 
S_3(x)&=&\frac{\pi^2}{6}x-\frac{\pi}{4}x^2+\frac{1}{12}x^3,\quad 0\le x\le 2\pi.
\\
C_4(x)&=&\frac{\pi^4}{90}-\frac{\pi^2}{12}x^2+\frac{\pi}{12}x^3-\frac{1}{48}x^4,\quad 0\le x\le 2\pi.
\end{eqnarray} 
\end{exa}

\subsection{Non-polynomial Branch}\label{sec.Cl}
The functions $\Cl_j(x)$ are represented by truncated Chebyshev
series as described in the following two subsections.
\subsubsection{$j=2$}
The function $\Cl_2(x)=S_2(x)$ is represented in power series as
\cite{AshourMComp10}
\begin{equation}
\Cl_2(x) = x-x\ln x +\frac{1}{2}x^3\sum_{n\ge 0} \frac{|B_{2n+2}|}{(n+1)(2n+3)!}x^{2n},\quad 0\le x\le \pi,
\label{eq.cl2tayl}
\end{equation}
and
\begin{equation}
\Cl_2(x) = (\pi-x)\left\{
\ln 2 -\sum_{n\ge 1} \frac{(2^{2n}-1)|B_{2n}|}{2n(2n+1)!}(\pi-x)^{2n}
\right\},\quad 0<x<2\pi.
\end{equation}
The current implementation follows precisely the separate Chebyshev expansions
in the subintervals $-\pi/2\le x \le \pi/2$ and $\pi/2 \le x \le 3\pi/2$
proposed by K\"olbig \cite{Kolbig}. The minor difference is that
the precision of the expansion coefficients has been increased to
33 digits with the aim to support quad precision calculations.
\begin{rem}
The GNU Scientific Library offers another implementation of $\Cl_2(x)$
by writing the non-logarithmic terms of (\ref{eq.cl2tayl})
as a product of $x$ times a series of 
even Chebyshev Polynomials in the interval $-\pi \le x \le  \pi$.
\end{rem}

\subsubsection{$j\ge 3$}
The functions $\Cl_j(x)$ for $3\le j \le 10$ have been bootstrapped
from a Chebyshev expansion of $\Cl_2(x)$ in the same spirit of 
repeated integration as for the polynomials in Section \ref{sec.poly} \cite{KalmykovCPC172}:
\begin{equation}
\Cl_{j+1}(x) = \Cl_{j+1}(0)-(-1)^j\int_0^x \Cl_{j}(x) dx.
\label{eq.iterI}
\end{equation}
Here $\Cl_j(0)=0$ if $j$ is even, $\Cl_j(0)=\zeta(j)$ if $j$ is odd.

To support a simple recurrence, the fundamental representation is a
Chebyshev expansion over an
entire period (similar to Wood's coverage \cite{WoodMathComp22}),
transforming (\ref{eq.cl2tayl}) into
\begin{equation}
\Cl_2(x) = -x \ln x +\sum_{r=1,3,5,7,\ldots} \alpha_r T_r(x/\pi),\quad
-\pi \le x \le \pi,
\label{eq.cl2cheb}
\end{equation}
with odd-indexed Chebyshev Polynomials $T$.
The $\alpha_r$ have been computed with Maple \cite{maple} by plugging
the formula which expands integer powers of $x$ in sums of $T_r(x)$
\cite{Cody,Fox,FraserJACM12}
into (\ref{eq.cl2tayl}):
\begin{equation}
\frac{\alpha_r}{\pi} = \delta_{1,r}
+\sum_{n= \max[0,\frac{r-3}{2}]}^\infty \binom{2n+3}{n+(3-r)/2}\frac{\zeta(2n+2)}{(n+1)(2n+3)2^{4n+4}},\, r=1,3,5,\ldots
\end{equation}
A table of explicit values of $\alpha_r$ is the
array \texttt{MooreT} in the file \texttt{anc/clausen.c} of
the ancillary material.

The leading logarithmic term in (\ref{eq.cl2cheb})
is repeatedly integrated via \cite[2.723]{GR}
\begin{equation}
\int_0^x [-\frac{1}{n!} \log(t)
+\alpha]t^n\,dt
=
\left[
-\frac{1}{(n+1)!} \log x
+\frac{1+\alpha (n+1)!}{(n+1)!(n+1)}
\right] x^{n+1}
.
\label{eq.logItr}
\end{equation}
\begin{exa}
\begin{eqnarray}
\int_0^x (-t)\ln t\,dt &=& -\frac{1}{2}x^2\ln x+\frac{x^2}{4},
\\
\int_0^x \left[-\frac{1}{2}t^2\ln t+\frac{t^2}{4}\right]dt
&=& -\frac{1}{6}x^3\ln x+\frac{5}{36}x^3,
\\
\int_0^x \left[ -\frac{1}{6}t^3\ln t+\frac{5}{36}t^3\right]dt &=& -\frac{1}{24}x^4\ln x +\frac{13}{288}x^4.
\end{eqnarray}
\end{exa}

The Chebyshev Polynomials in the integral kernel of the right
hand sides in (\ref{eq.iterI}) are mapped repeatedly onto other
Chebyshev Polynomials of adjacent indices with the generic
indefinite integrals \cite{ClenshawMPCPS53,PiessensJCAM121}
\begin{equation}
\int T_n(x) dx
=
\left\{
\begin{array}{ll}
\frac12\left[ \frac{T_{n+1}(x)}{n+1}-\frac{T_{|n-1|}(x)}{n-1}\right],& n=0\, \mathrm{or}\, n>1,\\
\frac12 \frac{T_{n+1}(x)}{n+1},& n=1.
\end{array}
\right.
\end{equation}
The associate definite integrals are
\begin{equation}
\int_0^x T_n(y) dy
=
\left\{
\begin{array}{ll}
 \frac{T_{n+1}(x)}{2(n+1)} -\frac{T_{|n-1|}(x)}{2(n-1)}, & n\, \mathrm{even}\, \ge 0, \\
 \frac{T_{n+1}(x)}{2(n+1)} -\frac{T_{n-1}(x)}{2(n-1)}+ (-)^{(n-1)/2}\frac{n}{n^2-1}, & n\, \mathrm{odd}\, \ge 3,\\
 \frac{T_{n+1}(x)}{2(n+1)} +\frac{T_{n-1}(x)}{4}, & n =1.\\
\end{array}
\right .
\label{eq.TItr}
\end{equation}

\section{Coding Details}
\subsection{Distribution}
The files \texttt{clausen.h}, \texttt{clausen.c} and \texttt{Makefile}
in the ancillary directory contain the source code. \texttt{clausen.h} is
the header file for inclusion in application programs; \texttt{clausen.c}
contains the three public functions and a test interface (included if
the preprocessor macro \texttt{TEST} is added), and \texttt{Makefile}
supports the primitive compiler call such that
\begin{verbatim}
make
\end{verbatim}
suffices to create the \texttt{clausen.o} linker file on standard Unices.
A C99 compiler is required because the program uses
\begin{enumerate}
\item
the \texttt{NAN} value to indicate out-of-range parameters,
\item
local variable definitions in loop constructs,
\item
and delayed variable declarations inside blocks.
\end{enumerate}

\subsection{clausen.c}
\subsubsection{clausen}
This function calculates of $\Cl_j(x)$
defined in (\ref{eq.Clausdef}) by dispatching it to calculations of (\ref{eq.Cdef})
or (\ref{eq.Sdef}) depending on the parity of $j$.

\subsubsection{clausenSin}
This function calculates $S_j(x)$
defined in (\ref{eq.Sdef}). $x$ is first reduced to the fundamental
period $0\le x< 2\pi$, then to the interval
$0\le x < \pi$ with the aid of $S_j(x)=-S_j(2\pi-x)$.
For $j=1$, (\ref{eq.S1}) is returned.
For large $j$, a partial sum of (\ref{eq.Sdef}) is accumulated with an upper limit
of the index $k$ fixed by a small lookup table depending on $j$.
The remaining cases are dispatched to \texttt{clausenPoly} if $j$ is odd
and to \texttt{clausenCheby} if $j$ is even.

\subsubsection{clausenCos}
This function calculates a value of $C_j(x)$
defined in (\ref{eq.Cdef}). $x$ is first reduced to the fundamental
period $0\le x< 2\pi$, then to the interval $0\le x< \pi$
via $C_j(x)=C_j(2\pi -x)$. For $j=1$, (\ref{eq.C1}) is returned.
If $x$ is a multiple of $\pi$---equivalent to an attempt to calculate $\zeta(1)$---\texttt{NAN} is returned
For large $j$, a partial sum of (\ref{eq.Cdef}) is accumulated with an upper limit
of the index $k$ fixed by a small lookup table depending on $j$.
The remaining cases are dispatched to \texttt{clausenPoly} if $j$ is even
and to \texttt{clausenCheby} if $j$ is odd.

\subsubsection{repeatLogInt}
This function computes the repeated integral
\begin{equation}
\int_0^x \int_0^{x'}\cdots \int_0^{x^{(m)}} [-x^{(m)}\ln x^{(m)}] dx^{(m)} \cdots dx' dx
\end{equation}
as required to lift the logarithmic term in (\ref{eq.cl2tayl}) from $\Cl_2$ to $\Cl_s$
with the aid of (\ref{eq.logItr}).

\subsubsection{clausenPoly}
This function calculates the sums satisfying the requirements
of Section \ref{sec.poly}. An initial polynomial with the coefficients
taken from (\ref{eq.C2pol}) is defined, and the chain of integrations 
rising the index is executed term by term. The $\zeta$-values are
read from a static table.

\subsubsection{evalCheby}
Given a vector of expansion coefficients, this function calculates
a Chebyshev series up to a maximum index implied by the length
of the vector at some fixed argument. Similar to the Horner scheme
of polynomial algebra, the terms at high indices are gathered first,
assuming that these carry the smallest numbers, and the
recurrence $T_{n+1}(x)=2xT_n(x)-T_{n-1}(x)$ is used to work backwards
through the vector of coefficients.

\subsubsection{clausenCheby}
This function covers the computation discussed in Section
\ref{sec.Cl}. The function $\Cl_2$ is handled exactly as
proposed by K\"olbig---this is not strictly needed, but
it is more efficient than using the generic evaluation because
the number of expansion coefficients is smaller
if $x$ is close to zero.

The main branch of the function for indices larger than two takes the vector
of the $\alpha_r$ in (\ref{eq.cl2cheb}) from a static list 
and inserts a zero at each second place representing the vanishing
coefficients at even indices. The step (\ref{eq.TItr})
is then repeatedly executed to lift the index from $\Cl_2\to \ldots \Cl_j \ldots \to \Cl_s$
one by one.
The sign in (\ref{eq.iterI}) is taken into account, and because
the arguments of the Chebyshev Polynomials are $x/\pi$,
a factor $\pi$ is emitted with each integration (implied by the substitution
$x/\pi \to x$ of the integration variable).
The constant of
integration is fixed by manipulation of the coefficient at index zero:
it either is zero if the intermediate $\Cl_j$
represents an even $j$, or $\zeta(j)$ if it represents an odd $j$.
Each integration is actually traced by modifying/mixing the vector of $\alpha_r$
in place. Once the index $s$ has been reached, the contribution 
from the intgrated $r$-sum of (\ref{eq.cl2cheb}) is gathered by evaluating
the final $\sum \alpha_rT_r$ by a call to \texttt{evalCheby}.

The contribution from the repeated integral of the logarithmic term
is added by calling \texttt{repeatLogInt}, noting again the sign switches
induced in front of the integral (\ref{eq.iterI}).

\bibliographystyle{amsplain}
\bibliography{all}

\end{document}